\documentclass[preprints,article,accept,moreauthors,pdftex]{my-mdpi}

\usepackage{epsfig,graphicx,epstopdf}

\firstpage{1} 
\pubvolume{xx}
\issuenum{1}
\articlenumber{5}
\pubyear{2019}
\copyrightyear{2019}

\externaleditor{Submitted: 31 Dec 2019; Revised: 21 and 28 Feb 2020; Accepted: 02 March 2020.}
\history{Published in: Mathematics 2020, Volume 8, Issue 3, 353.
DOI:{\changeurlcolor{black}\href{https://doi.org/10.3390/math8030353}{\detokenize{10.3390/math8030353}}}}

\pdfoutput=1

\usepackage{amsmath,amsfonts,amssymb}

\Title{The Stability and Stabilization of Infinite Dimensional Caputo-Time Fractional Differential Linear Systems}

\Author{Hanaa Zitane $^{1,\dagger}$\orcidA{}, 
Ali Boutoulout $^{2,\dagger}$\orcidB{} and 
Delfim F. M. Torres $^{3,\dagger}$*\orcidC{}}

\AuthorNames{Hanaa Zitane, Ali Boutoulout and Delfim F. M. Torres}

\address{%
$^{1}$ \quad MACS Laboratory, Department of Mathematics, Faculty of Sciences,\newline
University of Moulay Ismail, 11201 Meknes, Morocco; h.zitane@edu.umi.ac.ma\\
$^{2}$ \quad MACS Laboratory, Department of Mathematics, Faculty of Sciences,\newline
University of Moulay Ismail, 11201 Meknes, Morocco; a.boutoulout@fs-umi.ac.ma\\
$^{3}$ \quad Center for Research and Development in Mathematics and Applications (CIDMA),\newline
Department of Mathematics, University of Aveiro, 3810-193 Aveiro, Portugal; delfim@ua.pt}

\corres{Correspondence: delfim@ua.pt; Tel.: +351-234-370-668 (D.F.M.T.)}

\firstnote{The authors contributed equally to this work.} 
	
\abstract{We investigate the stability and stabilization concepts for 
infinite dimensional time fractional differential linear systems 
in Hilbert spaces with Caputo derivatives.
Firstly, based on a family of operators generated by strongly 
continuous semigroups and on a probability density function, 
we provide sufficient and necessary conditions for the 
exponential stability of the considered class of systems. 
Then, by assuming that the system dynamics is symmetric 
and uniformly elliptic and by using the properties of the 
Mittag--Leffler function, we provide sufficient conditions 
that ensure strong stability. Finally, we characterize 
an explicit feedback control that guarantees the strong stabilization 
of a controlled Caputo time fractional linear system 
through a decomposition approach. Some examples are 
presented that illustrate the effectiveness of our results.}



\begin{document}
	
% -------------------------------------------
	
\section{Introduction}

Fractional order calculus is a natural generalization of classical integer order calculus. 
It deals with integrals and derivatives of an arbitrary real or complex order.
Fractional order calculus has become very popular, in recent years, due to its 
demonstrated applications in many fields of applied sciences and engineering,  
such as the spread of contaminants in underground water, the charge transport 
in amorphous semiconductors, and diffusion of pollution in the atmosphere
\cite{der22,der2,Refff3a}. Because it generalizes and includes in the limit 
the integer order calculus, the fractional calculus has the potential to 
accomplish much more than what integer order calculus achieves \cite{Ref3a}.
In particular, it has proved to be a powerful tool to describe 
long-term memory and hereditary properties of various dynamical 
complex processes \cite{Reff3a}, diffusion processes, such as those found 
in batteries \cite{Refffff3a} and electrochemical and control processes \cite{Reffffff3a},
to model and control epidemics \cite{MR3999702,MR3980195} and mechanical properties 
of viscoelastic systems and damping materials, such as stress and strain \cite{Reffff3a}.

One can find in the literature several different fractional calculus.
Here we use the fractional calculus of Caputo, which was introduced 
by Michele Caputo in his 1967 paper \cite{MR2379269}. Such calculus
has appeared, in a natural way, for representing observed phenomena 
in laboratory experiments and field observations, where the mathematical 
theory was checked with experimental data. Indeed, the operator introduced 
by Caputo in 1967, and used by us in the present work, represents an observed linear 
dissipative mechanism phenomena with a time derivative of order 0.15 
entering the stress-strain relation \cite{MR2379269}. More recently,
a variational analysis with Caputo operators has been developed, 
which provides further mathematical substance to the use
of Caputo fractional operators \cite{MR2984893,MR3443073}.

In the analysis and design of control systems, the stability issue 
has always an important role \cite{MR3887838,MR3783261}. 
For a dynamical system, an equilibrium state is said to be stable 
if such system remains close to this state for small disturbances, 
and for an unstable system the question is how to stabilize it, 
especially by a feedback control law \cite{MR3220777}.
The stabilization concept for integer order systems and related problems 
has been considered in several works, see, e.g., \cite{Prit,Curtain2,Trigg,Balak} 
and references cited therein. In \cite {Prit}, the relationship between 
the asymptotic behavior of a system, the spectrum properties of its dynamics, 
and the existence of a Lyapunov functional is provided. Several techniques 
are considered to study different kinds of stabilization, 
for example, the exponential stabilization is studied via a decomposition method \cite {Trigg} 
while the strong stabilization is developed using the Riccati approach \cite {Balak}.

Similarly as classical dynamical systems, stability analysis is a central task 
in the study of fractional dynamical systems, which has attracted 
increasing interest of many researchers \cite{MR3980195,MR3789859}. 
For finite dimensional systems, the stability concept 
for fractional differential systems equipped with the Caputo derivative
is investigated in many works \cite{MR2821248}.
In \cite{Mati}, Matignon studies the asymptotic behavior 
for linear fractional differential systems with the Caputo derivative, 
where their dynamics $A$ is a constant coefficient matrix. In this case 
the stability is guaranteed if the eigenvalues of the dynamics matrix $A$, 
$\lambda \in \sigma(A)$, satisfy $|arg(\lambda)|>\dfrac{\pi \alpha}{2}$ \cite{Mati}. 
Since then, many scholars have carried out further studies on the stability 
for different classes of fractional linear systems \cite{Qian,Chen}. 
In \cite{Qian}, stability theorems for fractional differential systems, which include 
linear systems, time-delayed systems, and perturbed systems, are established,
while in \cite{Chen}, Ge, Chen and Kou provide results on the Mittag--Leffler stability 
and propose a Lyapunov direct method, which covers 
the power law stability and the exponential stability. 
See also \cite{MR3844433}, where the Mittag--Leffler and the class-K function stability 
of fractional differential equations of order $\alpha \in (1,2)$ are investigated.
In 2018, the notion of regional stability was introduced 
for fractional systems in \cite{Chen2}, where the authors study 
the Mittag--Leffler stability and the stabilization of systems with Caputo derivatives,
but only on a sub-region of its geometrical domain. 
More recently, fractional output stabilization problems for distributed systems 
in the Riemann--Liouville sense were studied \cite{zit1,zitH,MR3994712}, 
where feedback controls, which ensure exponential, strong, and weak stabilization 
of the state fractional spatial derivatives, with real and complex orders, 
are characterized. 

An analysis of the literature shows that existing results 
on stability of fractional systems are essentially limited 
to finite-dimensional fractional order linear systems, while results 
on infinite-dimensional spaces are a rarity.  
In contrast, here we investigate global stability 
and stabilization  of infinite dimensional fractional dynamical 
linear systems in the Hilbert space $L^{2}(\Omega)$
with Caputo derivatives of fractional order $0 < \alpha < 1$.
In particular, we characterize exponential and strong stability
for fractional Caputo systems on infinite-dimensional spaces.

The remainder of this paper is organized as follows. 
In Section~\ref{sec:2}, some basic knowledge of fractional 
calculus and some preliminary results, which will be used 
throughout the paper, are given. In Section~\ref{sec:3}, 
we prove results on the global asymptotic and exponential stability 
of Caputo-time fractional differential linear systems.
In contrast with available results in the literature,
which are restricted to systems of integer order or 
to fractional systems in the finite dimensional state 
space $\mathbb{R}^{n}$, here we study a completely different 
class of systems: we investigate fractional linear systems 
where the state space is the Hillbert space $L^{2}(\Omega)$.
We also characterize the stabilization of a controlled 
Caputo diffusion linear system via a decomposition method.   
Section~\ref{sec:4} presents the main conclusions of the work 
and some interesting open questions that deserve further investigations.

% -------------------------------------------
	
\section{Preliminaries and Notation}
\label{sec:2}

In this section, we introduce several definitions 
and results of fractional calculus
that are used in the sequel.	

\begin{Definition}[\cite{der2}]
Let $0<\alpha <1$ and $T>0$. The Caputo derivative of fractional order $\alpha$ 
for an absolutely continuous function $y(\cdot)$ on $[0, T]$ can be defined as follows:
$$
^{C}D_{t}^{\alpha}y(t)=\dfrac{1}{\Gamma(1-\alpha)}
\int_{0}^{t}(t-s)^{-\alpha} \dfrac{d}{ds}y(s)\, \mathrm{d}s,
$$
where $\Gamma(1-\alpha)$ is the Euler Gamma function.
\end{Definition}

\begin{Lemma}[\cite{sol1}]
\label{solmild}
For any given function $g \in L^{2}(0, T, L^{2}(\Omega))$, we say that function 
$y \in C(0, T, L^{2}(\Omega))$ is a mild solution of the system
\begin{equation}
\label{system000}
\left\{
\begin{array}{ll}
^{C}D_{t}^{\alpha}y(t)=Ay(t)+g(t) & t\in]0,+\infty[\\
y(0)=y_{0}  & y_{0} \in L^{2}(\Omega)
\end{array}
\right.
\end{equation}
if it satisfies 
\begin{equation}
\label{sys000}
y(t)=S_{\alpha}(t)y_{0}
+\int_{0}^{t}(t-s)^{\alpha-1}K_{\alpha}(t-s)g(s) \, \mathrm{d}s,
\end{equation}
where 
\begin{equation}
\label{semigalpha}
S_{\alpha}(t)=\int_{0}^{+\infty}
\Psi_{\alpha}(\theta)S(t^{\alpha}\theta) \, \mathrm{d}\theta
\end{equation}
and
\begin{equation}
\label{semigKalpha}
K_{\alpha} (t)=\alpha\int_{0}^{+\infty}\theta
\Psi_{\alpha}(\theta)S(t^{\alpha}\theta) \, \mathrm{d}\theta
\end{equation}
with
\begin{equation}
\label {thetaa}
\Psi_{\alpha}(\theta)=\frac{1}{\alpha}\theta^{-1-\frac{1}{\alpha}}
T_{\alpha}(\theta^{-\frac{1}{\alpha}}), 
\end{equation} 
$(S(t))_{t\geq0}$ the strongly continuous semigroup generated 
by operator $A$, and $T_{\alpha}$ the probability 
density function defined on $(0,\infty)$ by
$$
T_{\alpha}=\frac{1}{\pi}\sum_{n=1}^{+\infty}(-1)^{n}\theta^{\alpha n-1}
\frac{\Gamma(n\alpha +1)}{n!}\sin(n\pi \alpha). 
$$	
\end{Lemma}

\begin{Remark}[\cite{proba}]
The probability density function $T_{\alpha}$ defined on $(0,\infty)$ 
satisfies
$$
T_{\alpha}(\theta)\geq0, \quad \theta \in(0,\infty), 
\quad \text{ and }
\int_{0}^{+\infty}T_{\alpha}(\theta)\, \mathrm{d}\theta=1.
$$
\end{Remark}

\begin{Definition}[\cite{mettagdef}]
The Mittag--Leffler function of one parameter is defined as
$$
E_{\eta}(z)=\sum_{n=0}^{+\infty}\frac{z^{n}}{\Gamma(\eta n+1)}
\quad \text{with} \quad Re(\eta)>0, \quad z~\in \mathbb{C}.
$$
\end{Definition}

\begin{Definition}[\cite{mettagdef}]
The Mittag--Leffler function of two parameters is defined as
$$
E_{\eta,\beta}(z)
=\sum_{n=0}^{+\infty}\frac{z^{n}}{\Gamma(\eta n+ \beta)}
\quad \text{with} \quad  Re(\eta)>0, \quad \beta >0, 
\quad z \in \mathbb{C}.
$$	
\end{Definition}

\begin{Remark}
The Mittag--Leffler function appears naturally in the solution 
of fractional differential equations and in various applications: 
see \cite{mettagdef} and references therein. The exponential function 
is a special case of the Mittag--Leffler function \cite{MR3914412}:
for $\beta=1$ one has $E_{\eta,1}(z)=E_{\eta}(z)$ 
and $E_{1,1}(z)=e^{z}$.		
\end{Remark}

\begin{Lemma}[\cite{MAINARDI2}]
\label{mettaglem}
The Mittag--Leffler function $E_{\alpha}(-t^{\alpha})$ is
completely monotonic: for all $0 <\alpha < 1$, 
for all $n\in \mathbb{N}$ and $t > 0$, one has
$$
(-1)^{n}\frac{d^{n}}{dt^{n}} E_{\alpha}(-t^{\alpha})\geq 0.
$$
\end{Lemma}

\begin{Lemma}[\cite{mettaglemma2}]
\label{mettaglem2a}
The generalized Mittag--Leffler function $E_{\eta, \beta}(-x)$, $x \geq 0$, 
is completely monotonic for $\eta, \beta >0$  if and only if 
$\eta \in (0, 1]$ and $\beta \geq \eta$.	
\end{Lemma}

\begin{Lemma}[\cite{der1}]
\label{mettaglem2b}
Let $\beta >0$, $0<\eta < 2$, and $\mu$ be an arbitrary real number 
such that $\frac{\pi \eta}{2}< \mu <\min\{\pi, \pi \eta\}$. Then, 
the following asymptotic expressions hold: 
\begin{itemize}
\item if $|arg(z)|\leq \mu$ and $|z|>0$, then
\begin{equation}
|E_{\eta,\beta}(z)|\leq M_{1}(1+|z|)^{(1-\beta)/\eta}e^{Re(z^{\frac{1}{{\eta}}})}
+\frac{M_{2}}{1+ |z|};
\end{equation}
\item if $ \mu <|arg(z)|\leq \pi$ and $|z|\geq 0$, then
\begin{equation}
\label{theo2}
|E_{\eta,\beta}(z)|\leq \frac{M_{2}}{1+ |z|},
\end{equation}
\end{itemize}
where $M_{1}$ and $M_{2}$ are positive constants.
\end{Lemma} 

% --------------------------------------------------------

\section{Main Results}
\label{sec:3}

Our main goal is to study the stability and provide stabilization 
for a class of abstract Caputo-time fractional differential linear systems.

%--------------------------------------------------------

\subsection{Stability of Time Fractional Differential Systems}
\label{sec3.1}
 
Let $\Omega$ be an open bounded subset of $\mathbb{R}^{n}$, $n=1,2,3,\ldots$, 
and let us consider the following abstract time fractional order  differential system:
\begin{equation}
\label{system0}
\left\{
\begin{array}{ll}
^{C}D_{t}^{\alpha}z(t)=Az(t), & t \in \, ]0,+\infty[,\\
z(0)=z_{0},  & z_{0} \in L^{2}(\Omega),
\end{array}
\right.
\end{equation}
where $^{C}D_{t}^{\alpha}$ is the  left-sided Caputo fractional derivative of order 
$0 < \alpha < 1$, the second order operator 
$A: D(A) \subset L^{2}(\Omega) \longrightarrow L^{2}(\Omega)$ is
linear, with dense domain and such that the coefficients do not depend
on time $t$, and such that it is also the infinitesimal generator 
of the $C_{0}$-semi-group $(S(t))_{t \geq 0}$ on the Hilbert state space 
$L^{2}(\Omega)$ endowed with its usual inner product  $<\cdot,\cdot>$ 
and the corresponding norm $\Arrowvert \cdot \Arrowvert$. 
The unique mild solution of system \eqref{system0} can be written, 
from Lemma~\ref{solmild}, as 
$$
z(t)=S_{\alpha}(t)z_{0},
$$
where $S_{\alpha}(t)$ is defined by \eqref{semigalpha}.

We begin by proving the following lemma, which will be used thereafter.

\begin{Lemma}
\label{lemm1}
Let $A$ be the infinitesimal generator of a $C_{0}$-semi-group 
$(S(t))_{t \geq 0}$ on the Hilbert space $L^{2}(\Omega)$. 
Assume that there exists a function $h(\cdot) 
\in L^{2}(0,+\infty;\mathbb{R}^{+})$ satisfying
\begin{equation}
\label{eq1}
\Arrowvert S_{\alpha}(t+s)z \Arrowvert 
\leq h(t)\Arrowvert S_{\alpha}(s)z\Arrowvert, 
\quad \forall \, t, s \geq 0,
\quad \forall \, z \in L^{2}(\Omega).
\end{equation} 
Then the operators $(S_{\alpha}(t))_{t \geq 0}$ are uniformly bounded.	
\end{Lemma}

\begin{proof}
To prove that $(S_{\alpha}(t))_{t \geq 0}$ are bounded, we have to show that  
\begin{equation}
\label{Salpha}
\forall  z \in L^{2}(\Omega) ~\underset{t\geq 0}{\sup}
\Arrowvert S_{\alpha}(t)z \Arrowvert < \infty.
\end{equation}
By \emph{reductio ad absurdum}, let us suppose that \eqref{Salpha} does not hold, 
which means that there exists a sequence $(t_{s}+\tau_{n})$,  
$t_{s}> 0$ and $\tau_{n} \longrightarrow +\infty$, satisfying
\begin{equation}\label{contr}
\Arrowvert S_{\alpha}(t_{s}+\tau_{n})z \Arrowvert 
\longrightarrow +\infty ~as  ~n \longrightarrow +\infty.	
\end{equation}
From relation
$$ 
\displaystyle \int_{0}^{+\infty} \Arrowvert S_{\alpha}(s+\tau_{n})z\Arrowvert^{2} 
\, \mathrm{d}s= \displaystyle \int_{\tau_{n}}^{+\infty} \Arrowvert 
S_{\alpha}(s)z\Arrowvert^{2} \, \mathrm{d}s,~~0\leq s<+\infty,
$$
it follows that the right-hand side goes to $0$ as $n \longrightarrow +\infty$.
Using Fatou's Lemma yields  
$$
\underset{n \longrightarrow +\infty }{\lim ~ \inf}~\Arrowvert 
S_{\alpha}(s+\tau_{n})z\Arrowvert=0~~ \forall~ s>0.
$$
Therefore, for some $s_{0}<t_{s}$, we may find a subsequence 
$\tau_{n_{k}}$ such that $$\underset{k \longrightarrow +\infty }{\lim}
\Arrowvert S_{\alpha}(s_{0}+\tau_{n_{k}})z\Arrowvert=0.
$$
By virtue of condition \eqref{eq1}, one obtains 
$$ 
\Arrowvert S_{\alpha}(t_{s}+\tau_{n_{k}})z \Arrowvert 
\leq h(t_{s}-s_{0})\Arrowvert S_{\alpha}(s_{0}+\tau_{n_{k}})z\Arrowvert 
\underset{k \longrightarrow +\infty}{\longrightarrow 0},
$$
which contradicts \eqref{contr}. The intended conclusion
follows from the uniform boundedness principle.
\end{proof}

\begin{Definition}
Let $z_{0} \in L^{2}(\Omega)$. System \eqref{system0} is said 
to be exponentially stable if there exist two strictly positive
constants, $M>0$ and $\omega> 0$, such that
$$ 
\Arrowvert z(t) \Arrowvert \leq Me^{-\omega t} \Arrowvert z_{0} \Arrowvert, 
\quad \forall t \geq 0.
$$  
\end{Definition}

The next theorem provides necessary and sufficient conditions 
for exponential stability of the abstract fractional order  
differential system \eqref{system0}.

\begin{Theorem}
\label{prop3}  
Suppose that the operators $(S_{\alpha}(t))_{t \geq 0}$ 
fulfill assumption \eqref{eq1} and
\begin{equation}
\label{eq2}
\forall z \in L^{2}(\Omega) ~~~~~~  \Arrowvert S_{\alpha}(t+s)z \Arrowvert 
\leq \Arrowvert S_{\alpha}(t)z\Arrowvert \cdot \Arrowvert S_{\alpha}(s)z\Arrowvert, 
\quad \forall t,s \geq 0.
\end{equation} 
Then, system \eqref{system0} is exponentially stable if, and only if, 
for every $z \in L^{2}(\Omega)$ there exists a positive constant 
$\delta< \infty$ such that 
\begin{equation}
\label{Condexp}
\displaystyle \int_{0}^{+\infty} \Arrowvert 
S_{\alpha}(t)z\Arrowvert^{2} \, \mathrm{d}t <\delta. 	
\end{equation}
\end{Theorem}	

\begin{proof}
One has
\begin{equation*}
\begin{split}
t\Arrowvert S_{\alpha}(t)z\Arrowvert^{2} 
&= \displaystyle \int_{0}^{t} \Arrowvert S_{\alpha}(t)z\Arrowvert^{2} \, \mathrm{d}s\\
&= \displaystyle \int_{0}^{t} \Arrowvert S_{\alpha}(t-s+s)z\Arrowvert^{2} \, \mathrm{d}s.
\end{split}
\end{equation*}
Combining assumption \eqref{eq1}, Lemma~\ref{lemm1}, 
and condition \eqref{Condexp}, one gets
\begin{equation*}
\begin{split}
t\Arrowvert S_{\alpha}(t)z \Arrowvert^{2}
&\leq  \displaystyle \int_{0}^{t} 
\Arrowvert S_{\alpha}(s)z\Arrowvert^{2} \Arrowvert 
S_{\alpha}(t-s)z\Arrowvert^{2} \, \mathrm{d}s\\
&\leq N\delta \Arrowvert z\Arrowvert^{2}
\end{split}
\end{equation*}
for some $N>0$. Therefore, for $t$ sufficiently large, it follows that 
$$
\Arrowvert S_{\alpha}(t)\Arrowvert<1.
$$ 	
Then, there exists $t_{1}>0$ such that
$$
\ln\Arrowvert S_{\alpha}(t)\Arrowvert<0,
\quad \forall t\geq t_{1}. 
$$ 
Thus, 
$$
\omega_{0}=\underset{t\geq 0}{\inf}~\dfrac{\ln\Arrowvert S_{\alpha}(t)\Arrowvert}{t} < 0.
$$
Now, let us show that
\begin{equation}
\label{111}
\omega_{0}= \underset{t\longrightarrow +\infty }{\lim}
\frac{\ln\Arrowvert S_{\alpha}(t)\Arrowvert}{t}.
\end{equation}
Let $t_{s}>0$ be a fixed number and 
$N^{'}=\underset{t \in [0, t_{s}] }{\sup}{\Arrowvert S_{\alpha}(t)\Arrowvert}$. 
Thus, for each $t>t_{s}$, there exists $m \in \mathbb{N}$ such that 
$mt_{s}\leq t \leq (m+1)t_{s}$. From \eqref{eq2}, it follows that
\begin{equation*}
\begin{split}
\Arrowvert S_{\alpha}(t)\Arrowvert 
&= \Arrowvert S_{\alpha}(mt_{s}+(t-mt_{s}))\Arrowvert\\
&\leq \Arrowvert S_{\alpha}(mt_{s})\Arrowvert\Arrowvert 
S_{\alpha}(t-mt_{s})\Arrowvert,
\end{split}
\end{equation*}
which yields
$$
\frac{\ln\Arrowvert S_{\alpha}(t)\Arrowvert}{t}
\leq \frac{\ln\Arrowvert S_{\alpha}(mt_{s})\Arrowvert}{t} 
+\frac{\ln\Arrowvert S_{\alpha}(t-mt_{s})\Arrowvert}{t}.
$$
Using again \eqref{eq2}, it results that
$$
\frac{\ln\Arrowvert S_{\alpha}(t)\Arrowvert}{t}
\leq\frac{mt_{s}}{t} \frac{\ln\Arrowvert 
S_{\alpha}(t_{s})\Arrowvert}{t_{s}} 
+\frac{\ln\Arrowvert N^{'}\Arrowvert}{t}.
$$
Since $\dfrac{mt_{s}}{t}\leq 1$ and $t_{s}$ is arbitrary, one obtains  
$$
\underset{t\longrightarrow +\infty }{\lim~ \sup}~
\frac{\ln\Arrowvert S_{\alpha}(t)\Arrowvert}{t}
\leq \underset{t>0 }{\inf}~ \frac{\ln\Arrowvert S_{\alpha}(t)\Arrowvert}{t}
\leq \underset{t\longrightarrow +\infty }{\lim~ \inf}
~\frac{\ln\Arrowvert S_{\alpha}(t)\Arrowvert}{t}.
$$
Consequently, \eqref{111} holds. Hence, we conclude that for all 
$\omega \in \, ]0, -\omega_{0}[$, there exists $M>0$ such that
$$
\forall z \in L^{2}(\Omega) 
\quad \Arrowvert S_{\alpha}(t)z\Arrowvert 
\leq Me^{-\omega t}\Arrowvert z\Arrowvert,~~ \forall t\geq0,
$$
which means that system \eqref{system0} is exponentially stable.
The converse is obvious.
\end{proof}	

\begin{Remark}
When $\alpha=1$, the conditions \eqref{eq1} and \eqref{eq2} are verified, 
and we retrieve from our Theorem~\ref{prop3}
the results established in \cite{Curtain2,Prit} about the 
exponential stability of system \eqref{system0} on $\Omega$,
which is equivalent to 
$$ 
\int_{0}^{+\infty} \Arrowvert  S(t)z\Arrowvert^{2} \, \mathrm{d}t 
<\infty, \quad \forall z \in L^{2}(\Omega).
$$ 
\end{Remark}

 \begin{Definition}
 Let $z_{0} \in L^{2}(\Omega)$. System \eqref{system0} is said to be strongly stable
 if its corresponding solution $z(t)$ satisfies
 $$
 \Arrowvert z(t) \Arrowvert \longrightarrow 0~~  as~~  t \longrightarrow +\infty.
 $$ 
 \end{Definition}

In our next theorem, we provide sufficient conditions 
that guaranty the strong stability of the fractional 
order differential system \eqref{system0}. 
The result generalizes the asymptotic result established by Matignon 
for finite dimensional state spaces, where the dynamics of the system $A$ 
is considered to be a  matrix with constant coefficients in $\mathbb{R}^{n}$ \cite{Mati}. 
In contrast, here we tackle the stability for a different class of systems. 
Precisely, we consider fractional systems where the system dynamics $A$ 
is a linear operator generating a strongly continuous semigroup 
in the infinite dimensional state space $L^{2}(\Omega)$.

\begin{Theorem}
\label{theorem2} 
Let $(\lambda_{p})_{p\geq1}$ and $(\phi_{p})_{p\geq1}$ 
be the eigenvalues and the corresponding eigenfunctions 
of operator $A$ on $L^{2}(\Omega)$. If $A$ is a symmetric 
uniformly elliptic operator, then system \eqref{system0} 
is strongly stable on $\Omega$.
\end{Theorem}

\begin{proof}
Since $A$ is a symmetric uniformly elliptic operator, 
it follows that system \eqref{system0} admits a weak 
solution defined by 
$$ 
z(t)= \sum_{p=1}^{+\infty} E_{\alpha} (\lambda_{p} t^{\alpha}) 
\langle z_{0}, \phi_{p}\rangle\phi_{p}
\quad \forall \, z_{0}\in L^{2}(\Omega),
$$ 
where $(\lambda_p)_{p\geq1}$ satisfy
$$
0>\lambda_{1}\geq \lambda_{2}  \geq \cdots \geq \lambda_{j} 
\geq \cdots, \lim\limits_{j \longrightarrow \infty}=-\infty,
$$ 
and $(\phi_{p})_{p\geq1}$ forms an orthonormal basis 
in $L^{2}(\Omega)$ \cite{solmittag,Aelliptic}.
Using the fact that function $E_{\alpha}(-t^{\alpha})$ is
completely monotonic, for all $ \alpha \in (0, 1)$ and $t > 0$ 
(Lemma~\ref{mettaglem}), yields
\begin{equation*}
\begin{split}
\Arrowvert z(t) \Arrowvert
&= \left\Arrowvert \displaystyle\sum \limits_{{p=1}}^{+\infty} 
E_{\alpha} (\lambda_{p} t^{\alpha}) \langle z_{0},
\phi_{p}\rangle\phi_{p}\right\Arrowvert\\
&\leq |E_{\alpha}(\lambda_{1} t^{\alpha})| 
\Arrowvert z_{0} \Arrowvert.
\end{split}
\end{equation*}
Moreover, from Lemma~\ref{mettaglem2b}, it follows that 
$$ 
\Arrowvert z(t) \Arrowvert \leq \frac{M_{2}}{1-\lambda_{1}t^{\alpha}}
\Arrowvert z_{0} \Arrowvert \longrightarrow 0
\quad \text{ as } \quad  t \longrightarrow +\infty
$$
for some $M_{2}>0$. Hence, 
system \eqref{system0} is strongly stable on $\Omega$. 
\end{proof}

\begin{Example}	
Let us consider, on $\Omega=]0, 1[$, the following 
one-dimensional fractional diffusion system defined by
\begin{equation}
\label{system55}
\left\{
\begin{array}{lll}
^{C}D_{t}^{0.5}z(x, t)=\dfrac{\partial^2 z}{\partial x^2}(x, t),
& x \in \Omega, \quad t \in \, ]0,+\infty[,\\
z(0, t)=z(1, t)=0,  &\forall t>0,\\
z(x,0)=z_{0}, &x \in \Omega,
\end{array}
\right.
\end{equation}
where the second order operator $A=\dfrac{\partial^2}{\partial x^2}$ 
has its spectrum given by the eigenvalues $\lambda_{p}=-(p\pi)^{2}$, $p \geq 1$, 
and the corresponding eigenfunctions are $\phi_{p}(x)=\sqrt{\frac{2}{1+(p\pi)^{2}}}\sin(p\pi x)$, 
$p\geq1$. Operator $A$ generates a $C_{0}$-semi-group $(S(t))_{t\geq0}$ defined by
$$ 
S(t)z_{0}=\sum_{p=1}^{+\infty}e^{\lambda_{p} t}\langle z_{0}, \phi_{p}\rangle\phi_{p}. 
$$	
Moreover, the solution of system \eqref{system55} is given by 
$$
S_{0.5}(t)z_{0}=\sum_{p=1}^{+\infty} E_{0.5} (\lambda_{p} t^{0.5})
\langle z_{0}, \phi_{p}\rangle\phi_{p}.
$$	
One has that operator $A$ is symmetric and uniformly elliptic.
Consequently, from our Theorem~\ref{theorem2}, we deduce that 
system \eqref{system55} is strongly stable on $\Omega$.
This is illustrated numerically in Figure~\ref{fig:1} for 
$z(x,0) =  \sin(\pi x)$, $t = 0.1$, $t=0.15$, $t=0.2$, and $t=1.0$.
% --------------------------------
\begin{figure}[!t]
\centering
\includegraphics[scale=0.5]{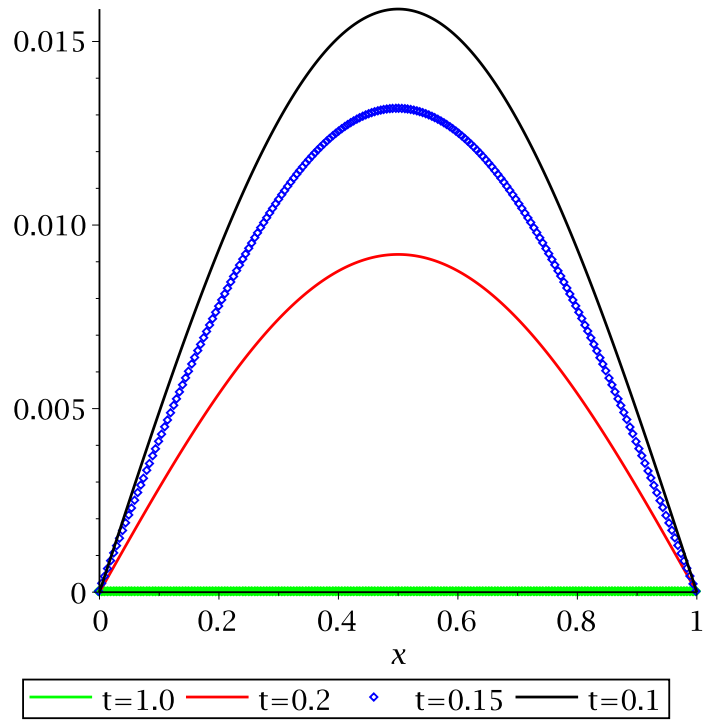}
\caption{The state of system \eqref{system55} for $z(x,0) =  \sin(\pi x)$,
$t = 0.1$, $t=0.15$, $t=0.2$, and $t=1.0$, illustrating the fact that 
\eqref{system55} is strongly stable on $\Omega=]0, 1[$.}\label{fig:1}
\end{figure}
% --------------------------------
\end{Example}

% --------------------------------------------------------

\subsection{Stabilization of Time Fractional Differential Systems}

Let $\Omega$ be an open bounded subset of $\mathbb{R}^{n}$, $n=1,2,3,\ldots$. 
We consider the following  Caputo-time fractional differential linear system:
\begin{equation}
\label{system0000}
\left\{
\begin{array}{ll}
^{C}D_{t}^{\alpha}z(t)=Az(t)+Bu(t), 
& t\in]0,+\infty[,\quad 0<\alpha<1,\\
z(0)=z_{0},  & z_{0} \in L^{2}(\Omega),
\end{array}
\right.
\end{equation}
with the same assumptions on $A$ as in Section~\ref{sec3.1}
and where $B$ is a bounded linear operator from $U$ into $L^{2}(\Omega)$, 
where $U$ is the space of controls, assumed to be a Hilbert space.
By Lemma~\ref{solmild}, the unique mild solution 
$z(\cdot)$ of system \eqref{system0000} is defined by 
\begin{equation}
\label{sys0000}
z(t)=S_{\alpha}(t)z_{0}+\displaystyle 
\int_{0}^{t}(t-s)^{\alpha-1}K_{\alpha}(t-s)Bu(s)\, \mathrm{d}s,
\end{equation}
where $S_{\alpha}(t)$ and $K_{\alpha}(t)$ are given, respectively, 
by \eqref{semigalpha} and \eqref{semigKalpha} .

\begin{Definition}
System \eqref{system0000} is said to be exponentially (respectively strongly) 
stabilizable if there exists a bounded operator $K \in \mathcal{L}(L^{2}(\Omega),U)$ 
such that the system	
\begin{equation}
\label{system2}
\left\{
\begin{array}{ll}
^{C}D_{t}^{\alpha}z(t)=(A+BK)z(t), &t\in \, ]0,+\infty[,\\
z(0)=z_0,   &z_0 \in L^{2}(\Omega),
\end{array}
\right.
\end{equation} 	
is exponentially (respectively strongly) stable on $\Omega$.
\end{Definition}

\begin{Remark}
It is clear that the exponential stabilization of system \eqref{system0000}  
implies the strong stabilization of \eqref{system0000}.
Note that the concept is general: when $\alpha=1$, we obtain 
the classical definitions of stability and stabilization.
\end{Remark}

Let $(S^{k}(t))_{t\geq0}$ be the strongly continuous semi-group generated by $A+BK$, 
where $K\in \mathcal{L}(L^{2}(\Omega),U)$ is the feedback operator. 
The unique mild solution of system \eqref{system0000} can be written as 
$$  
z(t)=S^{k}_{\alpha}(t)z_{0}    
$$
with 
$$
S^{k}_{\alpha}(t)= \displaystyle \int_{0}^{+\infty}
\Psi_{\alpha}(\theta)S^{k}(t^{\alpha}\theta) \, \mathrm{d}\theta,
$$
where $\Psi_{\alpha}(\theta)$ is defined by \eqref{thetaa}.

\begin{Theorem}
Let $A+BK$ generate a strongly continuous semi-group $(S^{k}(t))_{t \geq 0}$ on $L^{2}(\Omega)$. 
If the operator $(S^{k}_{\alpha}(t))_{t \geq 0}$ satisfies conditions \eqref{eq1} and \eqref{eq2} 
and if
$$ 
\forall z \in L^{2}(\Omega)~~~~~~\displaystyle \int_{0}^{+\infty} \Arrowvert 
S^{k}_{\alpha}(t)z\Arrowvert^{2} \, \mathrm{d}t <~ \infty 
$$	
holds, then system \eqref{system0000} is  exponentially stabilizable on $\Omega$.
\end{Theorem}	

\begin{proof}
The proof is similar to the proof of Theorem~\ref{prop3}.
\end{proof}	

\begin{Theorem}
\label{theo4} 
Let $(\lambda^{k}_{p})_{p\geq1}$ and $(\phi^{k}_{p})_{p\geq1}$ be the eigenvalues 
and the corresponding eigenfunctions of operator $A+BK$ on $L^{2}(\Omega)$.
If $A+BK$ is a symmetric uniformly elliptic operator, 
then system \eqref{system0000} is strongly stabilizable on $\Omega$.
\end{Theorem}

\begin{proof}
The proof is similar to the proof of Theorem~\ref{theorem2}.
\end{proof}	

\begin{Example}	
Let us consider, on $\Omega=]0, 1[$, 
the following fractional differential
system of order $\alpha=0.2$:
\begin{equation}
\label{system55b}
\begin{cases}
^{C}D_{t}^{0.2}z(x, t)=
\displaystyle \frac{1}{100}\dfrac{\partial^2 z}{\partial x^2}(x, t)
+\displaystyle \frac{1}{2} z(x,t)
+ B Kz(x, t), &(x,t) \in \Omega \times ]0,+\infty[,\\
z(0, t) = z(1, t) = 0,  &\forall t>0,\\
z(x,0)=z_{0}, &x \in \Omega,
\end{cases}
\end{equation}	
with the linear bounded operator $B=I$ and where we take $K=-B^{*}=-I$. 
The operator 
$$
A+BK=\frac{1}{100}\dfrac{\partial^2}{\partial x^2}-\frac{1}{2}, 
$$
with spectrum given by the eigenvalues 
$\lambda^{k}_{p}= -\frac{1}{2} -\frac{1}{100}(p\pi)^{2}$, $p \geq 1$, 
and the corresponding eigenfunctions $\phi^{k}_{p}(x)
=\sqrt{\frac{2}{1+(p\pi)^{2}}}\cos(p\pi x)$, $p\geq1$,	
generates a $C_{0}$-semi-group $(S^{k}(t))_{t\geq 0}$ defined by
$$ 
S^{k}(t)z_{0}=\displaystyle\sum_{p=1}^{+\infty}e^{\lambda^{k}_{p} t}\langle z_{0}, 
\phi^{k}_{p}\rangle\phi^{k}_{p}. 
$$
Furthermore, the solution of system \eqref{system55b} can be written as
$$
z(t)=S^{k}_{0.2}(t)z_{0}=\displaystyle\sum_{p=1}^{+\infty} 
E_{0.2} (\lambda^{k}_{p} t^{0.2})\langle z_{0}, \phi^{k}_{p}\rangle\phi^{k}_{p}.
$$
It is clear that $A+BK$ is a symmetric and uniformly elliptic operator.	
Hence, from Theorem~\ref{theo4}, we deduce that system 
\eqref{system55b} is strongly stabilizable on $\Omega$, i.e., the system
\begin{equation*}
\begin{cases}
^{C}D_{t}^{0.2}z(x, t)=
\displaystyle \frac{1}{100}\dfrac{\partial^2 z}{\partial x^2}(x, t)
+\displaystyle \frac{1}{2} z(x,t) + B u(t), 
&(x,t) \in \Omega \times ]0,+\infty[,\\
z(0, t) = z(1, t) = 0,  &\forall t>0,\\
z(x,0)=z_{0}, &x \in \Omega,
\end{cases}
\end{equation*} 
is strongly stabilizable by the feedback control $u(t)=-B^{*}z(t)$.
Figure~\ref{fig:2} shows, for $z(x,0) = x (x-1)$, 
that the state $z(x,t)$ of system \eqref{system55b} 
is unstable at $t=0$. Moreover, we see that the state evolves close to 0 at $t=10$. 
Numerically, the state is stabilized by $u(t)=-B^{*}z(t)$
with an error equal to $1.75\times10^{-04}$.
% --------------------------------
\begin{figure}[!t]
\centering
\includegraphics[scale=0.7]{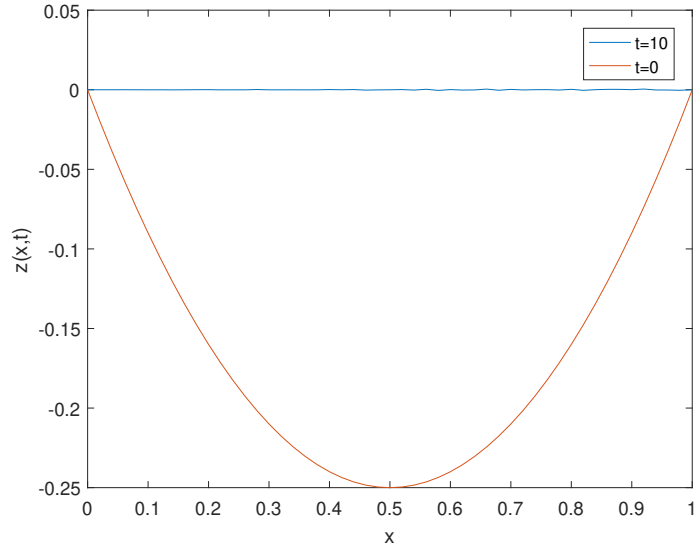}
\caption{The state of system \eqref{system55b} for 
$z(x, 0)=x(x-1)$, $t=0$, and $t=10$, illustrating the 
fact that \eqref{system55b} is unstable at $t=0$ 
but it is stabilized at $t=10$ on $\Omega=]0, 1[$.}\label{fig:2}
\end{figure}
% --------------------------------
\end{Example}

% --------------------------------------------------------

\subsection{Decomposition Method}

Now, we study the stabilization of system \eqref{system0000} 
using the decomposition method, which consists in decomposing 
the state space and the system using the spectral properties of operator $A$.

Let $\xi> 0$ be fixed and assume that there are at most finitely many 
nonnegative eigenvalues of $A$ and each with finite dimensional eigenspace. 
In other words, assume there exists $l \in \mathbb{N}$ such that 
\begin{equation}
\label{spectrumdecom}
\sigma(A)= \sigma_{u}(A)\cup \sigma_{s}(A),
\end{equation}
where $\sigma_{u}(A)= \sigma(A) \cap \{ \lambda_{p}, ~~~p=1, 2,\ldots,l \}$, 
$\sigma_{s}(A)= \sigma(A) \cap \{\lambda_{p}, ~~~ p=l+1, l+2 \ldots \}$ 
with $\lambda_{l} \geq 0$ and $\lambda_{l+1} \leq -\xi$. 
Because the sequence $(\phi_{p})_{p\geq1}$ forms a complete and orthonormal basis 
in $H=L^{2}(\Omega)$, it follows that the state space $H$ can be decomposed as
\begin{equation}
\label{eqa}
H=H_{u} \oplus H_{s},
\end{equation}
where $H_{u}=PH=\mathrm{span}\{ \phi_{1}, \phi_{2}, \ldots, \phi_{l} \}$ 
and $H_{s}=(I-P)H=\mathrm{span}\{ \phi_{l+1}, \phi_{l+2}, \ldots \}$ 
with $P \in \mathcal{L}(H)$ the projection operator \cite{decomp}.
Hence, system \eqref{system0000} can be decomposed 
into the following two sub-systems: 
\begin{equation}
\label{system3}
\begin{cases}
^{C}D_{t}^{\alpha}z_{u}(t)=A_{u}z_{u}(t)+PBu(t), \\
z_{0u}=Pz_0,
\end{cases}
\end{equation}
and 
\begin{equation}
\label{system4}
\begin{cases}^{C}D_{t}^{\alpha}z_{s}(t)
=A_{s}z_{s}(t)+(I-P)Bu(t), \\
z_{0s}=(I-P)z_0,
\end{cases}
\end{equation} 
where $A_{s}$ and $A_{u}$ are the restrictions of $A$ on $H_{s}$ 
and $H_{u}$, respectively, and are such that 
$\sigma(A_{s})=\sigma_{s}(A)$, $\sigma(A_{u})=\sigma_{u}(A)$, 
and $A_{u}$ is a bounded operator on $H_{u}$.

Our next result asserts that stabilization of system
\eqref{system0000} is equivalent to the one of system \eqref{system3}.

\begin{Theorem}
Let the spectrum $\sigma (A)$ of $A$ satisfy the above spectrum 
decomposition assumptions \eqref{spectrumdecom} for some $\xi> 0$ 
and $A_{s}$ be a symmetric uniformly elliptic operator. 
If system \eqref{system3} is strongly stabilizable by the control
\begin{equation}
\label{contDecom}
u(t)=D_{u}z_{u}(t)
\end{equation} 
with $D_{u}\in \mathcal{L}(H,U)$ such that
\begin{equation}
\label{stro22}
\Arrowvert z_{u}(t)\Arrowvert \leq C~t^{-\mu},
\quad \mu,C>0,	
\end{equation}	
then system \eqref{system0000} is strongly stabilizable 
using the feedback control $v(t)=D_{u}z_{u}(t)$.
\end{Theorem}	

\begin{proof}
Using the fact that system \eqref{system3} is strongly stabilizable 
by control \eqref{contDecom}, and inequality \eqref{stro22} yields
\begin{equation}
\label{infinity1}
\Arrowvert z_{u}(t)\Arrowvert \longrightarrow 0  \text{ as }  
t \longrightarrow +\infty
\end{equation}
and 
\begin{equation}
\label{conts}
\Arrowvert u(t)\Arrowvert 
\leq C \Arrowvert D_{u} \Arrowvert t^{-\mu},
\end{equation}
the unique weak solution of system \eqref{system4} can be written in the space $H_{s}$ as 	
$$ 
z_{s}(t)= \sum_{p=l+1}^{+\infty} E_{\alpha} (\lambda_{p} t^{\alpha}) 
\langle z_{0s}, \phi_{p}\rangle\phi_{p}+ \sum_{p=l+1}^{+\infty} 
\displaystyle \int_{0}^{t}(t-s)^{\alpha-1}
E_{\alpha,\alpha}(\lambda_{p}(t-s)^{\alpha})
\langle (I-P)Bu(s), \phi_{p}\rangle\phi_{p} \, \mathrm{d}s 
$$	
since $A_{s}$ is a symmetric uniformly elliptic operator \cite{solmittag}. 
Using the spectrum decomposition relation \eqref{spectrumdecom}, 
Lemma~\ref{mettaglem}, and Lemma~\ref{mettaglem2a}, one has that
\begin{equation}
\label{RA}
E_{\alpha} (\lambda_{p} t^{\alpha})\leq E_{\alpha}(- \xi t^{\alpha}) 
\quad \text{ for all } p\geq l+1
\end{equation}
and
\begin{equation}
\label{RA2}
E_{\alpha, \alpha}(\lambda_{p} (t-s)^{\alpha})
\leq E_{\alpha, \alpha}(- \xi (t-s)^{\alpha}) 
\quad \text{ for all } p\geq l+1.
\end{equation}
Then, feeding system \eqref{system4} by the same control 
$u(t)=D_{u}z_{u}(t)$ and using \eqref{conts}--\eqref{RA2}, 
it follows that 
\begin{equation*}
\begin{split}
\Arrowvert z_{s}(t) \Arrowvert
&\leq E_{\alpha}(- \xi t^{\alpha}) \Arrowvert z_{0s} \Arrowvert 
+ C \Arrowvert D_{u} \Arrowvert  \Arrowvert I-P \Arrowvert  \Arrowvert B \Arrowvert 
\displaystyle \int_{0}^{t}(t-s)^{\alpha-1}s^{-\mu}E_{\alpha,\alpha}(-\xi(t-s)^{\alpha}) \, \mathrm{d}s\\
&\leq E_{\alpha}(- \xi t^{\alpha}) \Arrowvert z_{0s} \Arrowvert 
+ C \Arrowvert D_{u} \Arrowvert  \Arrowvert I-P \Arrowvert  \Arrowvert B \Arrowvert
\displaystyle\sum \limits_{{n=1}}^{+\infty}\displaystyle\int_{0}^{t} 
\dfrac{(-\xi)^{n}(t-s)^{\alpha n+\alpha-1}s^{-\mu} \, \mathrm{d}s}{\Gamma(\alpha n+ \alpha)}\\
&\leq E_{\alpha}(- \xi t^{\alpha}) \Arrowvert z_{0s} \Arrowvert 
+ C \Arrowvert D_{u} \Arrowvert  \Arrowvert I-P \Arrowvert  \Arrowvert B \Arrowvert  
\displaystyle\sum \limits_{{n=1}}^{+\infty} 
\dfrac{(-\xi)^{n}t^{\alpha n+\alpha-\mu}}{\Gamma(\alpha n+\alpha-\mu-1)\Gamma(1-\mu)^{-1}}\\
&\leq E_{\alpha}(- \xi t^{\alpha}) \Arrowvert z_{0s} \Arrowvert 
+ C \Gamma(1-\mu)\Arrowvert D_{u} \Arrowvert  \Arrowvert I-P \Arrowvert  
\Arrowvert B \Arrowvert  t^{\alpha-\mu}
E_{\alpha,\alpha-\mu+1}(- \xi t^{\alpha}).
\end{split}
\end{equation*}
Lemma~\ref{mettaglem2b} implies that 
$$ 
\Arrowvert z_{s}(t) \Arrowvert \leq  \frac{M_{1}}{1+ \xi t^{\alpha}}
\Arrowvert z_{s0}\Arrowvert+ C \Gamma(1-\mu)\Arrowvert D_{u} \Arrowvert  
\Arrowvert I-P \Arrowvert  \Arrowvert B \Arrowvert 
\frac{M_{2}t^{\alpha-\mu}}{1+ \xi t^{\alpha}} 
$$
for some $M_{1}, M_{2}>0$. Therefore, 
\begin{equation}
\label{infinity2}
\Arrowvert z_{s}(t) \Arrowvert \longrightarrow 0
\text{ as } t \longrightarrow +\infty.
\end{equation}
On the other hand, we have that 
\begin{equation}
\label{infinity3}
\Arrowvert z(t) \Arrowvert= \Arrowvert z_{s}(t) + z_{u}(t) \Arrowvert 
\leq \Arrowvert z_{s}(t) \Arrowvert + \Arrowvert  z_{u}(t) \Arrowvert. 
\end{equation}
Combining \eqref{infinity1}, \eqref{infinity2}, and \eqref{infinity3}, 
we deduce the strong stabilization of system \eqref{system0000}.
\end{proof}  

% -------------------------------
  	 
\section{Conclusions and Future Work}
\label{sec:4}

We investigated the stability problem of infinite dimensional
time fractional differential linear systems under Caputo 
derivatives of order $\alpha \in (0,1)$, where the state space 
is the Hillbert space $L^{2}(\Omega)$. We proved necessary 
and sufficient conditions for exponential stability and 
obtained a characterization for the asymptotic stability, 
which is guaranteed if the system dynamics is symmetric 
and uniformly elliptic. Moreover, some stabilization 
criteria were also proved. Finally, we investigated 
the strong stabilization of the system via a decomposition 
method, where an explicit feedback control is obtained.  
Illustrative examples were given, showing the effectiveness 
of the theoretical results.
As future work, we intend to extend our work 
to the class of infinite dimensional time fractional differential 
nonlinear systems. Various other questions are still open and deserve 
further investigations, such as, studying boundary stability and gradient 
stability for time fractional differential linear systems
or considering the more recent notion of $\Lambda$-fractional derivative
\cite{ref01}, and thus obtaining a geometrical interpretation.

% --------------------------------------

\authorcontributions{Each author equally contributed to this paper, 
read and approved the final manuscript.}

% --------------------------------------

\funding{This research was funded by Moulay Ismail University (H.Z.);
by Hassan II Academy of Science and Technology, project N 630/2016 (A.B.);
and by The Portuguese Foundation for Science and Technology, R\&D unit CIDMA,
within project UIDB/04106/2020 (D.F.M.T.).}

% --------------------------------------

\acknowledgments{This research is part of first author's Ph.D. project,
which is carried out at Moulay Ismail University, Meknes, 
and has began during a one-month visit of Zitane 
to the R\&D Unit CIDMA, Department of Mathematics, 
University of Aveiro, Portugal, June 2019.
The hospitality of the host institution 
is here gratefully acknowledged.
The authors are strongly grateful to three anonymous referees 
for their suggestions and invaluable comments.}

% --------------------------------------

\conflictsofinterest{The authors declare no conflict of interest.} 

% --------------------------------------

\reftitle{References}

% --------------------------------------


\begin{thebibliography}{-------}
\providecommand{\natexlab}[1]{#1}
	
\bibitem[Rahimy(2010)]{der22}
Rahimy, M.
\newblock Applications of fractional differential equations.
\newblock {\em Appl. Math. Sci. (Ruse)} {\bf 2010}, {\em 4},~2453--2461.
	
\bibitem[Kilbas \em{et~al.}(2006)Kilbas, Srivastava, and Trujillo]{der2}
Kilbas, A.A.; Srivastava, H.M.; Trujillo, J.J.
\newblock {\em Theory and applications of fractional differential equations};
Elsevier Science B.V., Amsterdam,  2006.
	
\bibitem[Diethelm(2010)]{Refff3a}
Diethelm, K.
\newblock {\em The analysis of fractional differential equations}; Vol. 2004,
{\em Lecture Notes in Mathematics}, Springer-Verlag, Berlin,  2010.
\newblock
doi:{\changeurlcolor{black}\href{https://doi.org/10.1007/978-3-642-14574-2}{\detokenize{10.1007/978-3-642-14574-2}}}.
	
\bibitem[Hilfer(2000)]{Ref3a}
Hilfer, R.
\newblock {\em Applications of fractional calculus in physics}; World
Scientific Publishing Co., Inc., River Edge, NJ,  2000.
\newblock
doi:{\changeurlcolor{black}\href{https://doi.org/10.1142/9789812817747}{\detokenize{10.1142/9789812817747}}}.
	
\bibitem[Sabatier \em{et~al.}(2007)Sabatier, Agrawal, and Machado]{Reff3a}
Sabatier, J.; Agrawal, O.P.; Machado, J.A.T.
\newblock {\em Advances in fractional calculus}; Springer, Dordrecht,  2007.
\newblock
doi:{\changeurlcolor{black}\href{https://doi.org/10.1007/978-1-4020-6042-7}{\detokenize{10.1007/978-1-4020-6042-7}}}.
	
\bibitem[Gabano and Poinot(2011)]{Refffff3a}
Gabano, J.D.; Poinot, T.
\newblock Fractional modelling and identification of thermal systems.
\newblock {\em Signal Process.} {\bf 2011}, {\em 91},~531--541.
	
\bibitem[Ichise \em{et~al.}(1971)Ichise, Nagayanagi, and Kojima]{Reffffff3a}
Ichise, M.; Nagayanagi, Y.; Kojima, T.
\newblock An analog simulation of non-integer order transfer functions for
analysis of electrode processes.
\newblock {\em J. Electroanal. Chem. Interfacial Electrochem.} {\bf 1971}, {\em 33},~253--265.
	
\bibitem[Rosa and Torres(2019)]{MR3999702}
Rosa, S.; Torres, D.F.M.
\newblock Optimal control and sensitivity analysis of a fractional order {TB} model.
\newblock {\em Stat. Optim. Inf. Comput.} {\bf 2019}, {\em 7},~617--625.
\newblock
doi:{\changeurlcolor{black}\href{https://doi.org/10.19139/soic.v7i3.836}{\detokenize{10.19139/soic.v7i3.836}}}.
\newblock {\tt arXiv:1812.04507}

\bibitem[Silva and Torres(2019)]{MR3980195}
Silva, C.J.; Torres, D.F.M.
\newblock Stability of a fractional {HIV}/{AIDS} model.
\newblock {\em Math. Comput. Simulation} {\bf 2019}, {\em 164},~180--190.
\newblock
doi:{\changeurlcolor{black}\href{https://doi.org/10.1016/j.matcom.2019.03.016}{\detokenize{10.1016/j.matcom.2019.03.016}}}.
\newblock {\tt arXiv:1903.02534}
	
\bibitem[Bagley and Calico(1991)]{Reffff3a}
Bagley, R.L.; Calico, R.A.
\newblock Fractional order state equations for the control of viscoelastically damped structures.
\newblock {\em J. Guid. Control Dyn.} {\bf 1991}, {\em 14},~304--311.
	
\bibitem[Caputo(2008)]{MR2379269}
Caputo, M.
\newblock Linear models of dissipation whose {$Q$} is almost frequency independent. {II}.
\newblock {\em Fract. Calc. Appl. Anal.} {\bf 2008}, {\em 11},~4--14.
\newblock Reprinted from Geophys. J. R. Astr. Soc. {{\bf{1}}3} (1967), no.~5, 529--539.
	
\bibitem[Malinowska and Torres(2012)]{MR2984893}
Malinowska, A.B.; Torres, D.F.M.
\newblock {\em Introduction to the fractional calculus of variations}; Imperial
College Press, London,  2012; pp. xvi+275.
\newblock
doi:{\changeurlcolor{black}\href{https://doi.org/10.1142/p871}{\detokenize{10.1142/p871}}}.
	
\bibitem[Almeida \em{et~al.}(2015)Almeida, Pooseh, and Torres]{MR3443073}
Almeida, R.; Pooseh, S.; Torres, D.F.M.
\newblock {\em Computational methods in the fractional calculus of variations};
Imperial College Press, London,  2015; pp. xii+266.
\newblock
doi:{\changeurlcolor{black}\href{https://doi.org/10.1142/p991}{\detokenize{10.1142/p991}}}.
	
\bibitem[Mahmoud and Karaki(2018)]{MR3887838}
Mahmoud, M.S.; Karaki, B.J.
\newblock Improved stability analysis and control design of reset systems.
\newblock {\em IET Control Theory Appl.} {\bf 2018}, {\em 12},~2328--2336.
\newblock
doi:{\changeurlcolor{black}\href{https://doi.org/10.1049/iet-cta.2018.5410}{\detokenize{10.1049/iet-cta.2018.5410}}}.
	
\bibitem[Rocha \em{et~al.}(2018)Rocha, Silva, and Torres]{MR3783261}
Rocha, D.; Silva, C.J.; Torres, D.F.M.
\newblock Stability and optimal control of a delayed {HIV} model.
\newblock {\em Math. Methods Appl. Sci.} {\bf 2018}, {\em 41},~2251--2260.
\newblock
doi:{\changeurlcolor{black}\href{https://doi.org/10.1002/mma.4207}{\detokenize{10.1002/mma.4207}}}.
\newblock {\tt arXiv:1609.07654}
	
\bibitem[Sontag(2012)]{MR3220777}
Sontag, E.D.
\newblock Stability and feedback stabilization. In {\em Mathematics of
complexity and dynamical systems. {V}ols. 1--3}; Springer, New York,  2012; pp.~1639--1652.
\newblock
doi:{\changeurlcolor{black}\href{https://doi.org/10.1007/978-1-4614-1806-1_105}{\detokenize{10.1007/978-1-4614-1806-1_105}}}.
	
\bibitem[Pritchard and Zabczyk(1981)]{Prit}
Pritchard, A.J.; Zabczyk, J.
\newblock Stability and stabilizability of infinite-dimensional systems.
\newblock {\em SIAM Rev.} {\bf 1981}, {\em 23},~25--52.
\newblock
doi:{\changeurlcolor{black}\href{https://doi.org/10.1137/1023003}{\detokenize{10.1137/1023003}}}.
	
\bibitem[Curtain and Zwart(1995)]{Curtain2}
Curtain, R.F.; Zwart, H.
\newblock {\em An introduction to infinite-dimensional linear systems theory};
Vol.~21, {\em Texts in Applied Mathematics}, Springer-Verlag, New York, 1995.
\newblock
doi:{\changeurlcolor{black}\href{https://doi.org/10.1007/978-1-4612-4224-6}{\detokenize{10.1007/978-1-4612-4224-6}}}.
	
\bibitem[Triggiani(1975)]{Trigg}
Triggiani, R.
\newblock On the stabilizability problem in {B}anach space.
\newblock {\em J. Math. Anal. Appl.} {\bf 1975}, {\em 52},~383--403.
\newblock
doi:{\changeurlcolor{black}\href{https://doi.org/10.1016/0022-247X(75)90067-0}{\detokenize{10.1016/0022-247X(75)90067-0}}}.
	
\bibitem[Balakrishnan(1981)]{Balak}
Balakrishnan, A.V.
\newblock Strong stabilizability and the steady state {R}iccati equation.
\newblock {\em Appl. Math. Optim.} {\bf 1981}, {\em 7},~335--345.
\newblock
doi:{\changeurlcolor{black}\href{https://doi.org/10.1007/BF01442125}{\detokenize{10.1007/BF01442125}}}.
	
\bibitem[Wojtak \em{et~al.}(2018)Wojtak, Silva, and Torres]{MR3789859}
Wojtak, W.; Silva, C.J.; Torres, D.F.M.
\newblock Uniform asymptotic stability of a fractional tuberculosis model.
\newblock {\em Math. Model. Nat. Phenom.} {\bf 2018}, {\em 13},~Art. 9, 10.
\newblock
doi:{\changeurlcolor{black}\href{https://doi.org/10.1051/mmnp/2018015}{\detokenize{10.1051/mmnp/2018015}}}.
\newblock {\tt arXiv:1801.07059}
	
\bibitem[Zhang \em{et~al.}(2011)Zhang, Li, and Chen]{MR2821248}
Zhang, F.; Li, C.; Chen, Y.
\newblock Asymptotical stability of nonlinear fractional differential system
with {C}aputo derivative.
\newblock {\em Int. J. Differ. Equ.} {\bf 2011}, pp. Art. ID 635165, 12.
\newblock
doi:{\changeurlcolor{black}\href{https://doi.org/10.1155/2011/635165}{\detokenize{10.1155/2011/635165}}}.
	
\bibitem[Matignon(1996)]{Mati}
Matignon, D.
\newblock Stability results for fractional differential equations with
applications to control processing.
\newblock  Computational engineering in systems applications. Lille, France,
1996, Vol.~2, pp. 963--968.
	
\bibitem[Qian \em{et~al.}(2010)Qian, Li, Agarwal, and Wong]{Qian}
Qian, D.; Li, C.; Agarwal, R.P.; Wong, P.J.Y.
\newblock Stability analysis of fractional differential system with
{R}iemann-{L}iouville derivative.
\newblock {\em Math. Comput. Modelling} {\bf 2010}, {\em 52},~862--874.
\newblock
doi:{\changeurlcolor{black}\href{https://doi.org/10.1016/j.mcm.2010.05.016}{\detokenize{10.1016/j.mcm.2010.05.016}}}.
	
\bibitem[Li \em{et~al.}(2010)Li, Chen, and Podlubny]{Chen}
Li, Y.; Chen, Y.; Podlubny, I.
\newblock Stability of fractional-order nonlinear dynamic systems: {L}yapunov
direct method and generalized {M}ittag-{L}effler stability.
\newblock {\em Comput. Math. Appl.} {\bf 2010}, {\em 59},~1810--1821.
\newblock
doi:{\changeurlcolor{black}\href{https://doi.org/10.1016/j.camwa.2009.08.019}{\detokenize{10.1016/j.camwa.2009.08.019}}}.
	
\bibitem[Matar and Abu~Skhail(2018)]{MR3844433}
Matar, M.M.; Abu~Skhail, E.S.
\newblock On stability of nonautonomous perturbed semilinear fractional
differential systems of order {$\alpha\in(1,2)$}.
\newblock {\em J. Math.} {\bf 2018}, pp. Art. ID 1723481, 10.
\newblock
doi:{\changeurlcolor{black}\href{https://doi.org/10.1155/2018/1723481}{\detokenize{10.1155/2018/1723481}}}.
	
\bibitem[Ge \em{et~al.}(2018)Ge, Chen, and Kou]{Chen2}
Ge, F.; Chen, Y.; Kou, C.
\newblock {\em Regional analysis of time-fractional diffusion processes};
Springer, Cham,  2018.
\newblock
doi:{\changeurlcolor{black}\href{https://doi.org/10.1007/978-3-319-72896-4}{\detokenize{10.1007/978-3-319-72896-4}}}.
	
\bibitem[Zitane \em{et~al.}(2020)Zitane, Larhrissi, and Boutoulout]{zit1}
Zitane, H.; Larhrissi, R.; Boutoulout, A.
\newblock On the fractional output stabilization for a class of infinite
dimensional linear systems. In {\em Recent Advances in Modeling, Analysis and
Systems Control: Theoretical Aspects and Applications}; Springer, Cham, 2020; pp.~241--259.
	
\bibitem[Zitane \em{et~al.}(in press)Zitane, Larhrissi, and Boutoulout]{zitH}
Zitane, H.; Larhrissi, R.; Boutoulout, A.
\newblock Fractional output stabilization for a class of bilinear distributed systems.
\newblock {\em Rend. Circ. Mat. Palermo (2)} {\bf in press}.
\newblock
doi:{\changeurlcolor{black}\href{https://doi.org/10.1007/s12215-019-00429-w}{\detokenize{10.1007/s12215-019-00429-w}}}.
	
\bibitem[Zitane \em{et~al.}(2019)Zitane, Larhrissi, and Boutoulout]{MR3994712}
Zitane, H.; Larhrissi, R.; Boutoulout, A.
\newblock Riemann {L}iouville fractional spatial derivative stabilization of
bilinear distributed systems.
\newblock {\em J. Appl. Nonlinear Dyn.} {\bf 2019}, {\em 8},~447--461.
\newblock
doi:{\changeurlcolor{black}\href{https://doi.org/10.5890/JAND.2019.09.008}{\detokenize{10.5890/JAND.2019.09.008}}}.
	
\bibitem[Zhou and Jiao(2010)]{sol1}
Zhou, Y.; Jiao, F.
\newblock Existence of mild solutions for fractional neutral evolution equations.
\newblock {\em Comput. Math. Appl.} {\bf 2010}, {\em 59},~1063--1077.
\newblock
doi:{\changeurlcolor{black}\href{https://doi.org/10.1016/j.camwa.2009.06.026}{\detokenize{10.1016/j.camwa.2009.06.026}}}.
	
\bibitem[Mainardi \em{et~al.}(2007)Mainardi, Paradisi, and Gorenflo]{proba}
Mainardi, F.; Paradisi, P.; Gorenflo, R.
\newblock Probability distributions generated by fractional diffusion equations.
\newblock arXiv preprint {\tt arXiv:0704.0320} {\bf 2007}.
	
\bibitem[Erd\'{e}lyi \em{et~al.}(1981)Erd\'{e}lyi, Magnus, Oberhettinger, and Tricomi]{mettagdef}
Erd\'{e}lyi, A.; Magnus, W.; Oberhettinger, F.; Tricomi, F.G.
\newblock {\em Higher transcendental functions. {V}ol. {III}}; Robert E.
Krieger Publishing Co., Inc., Melbourne, Fla.,  1981.
	
\bibitem[Joshi \em{et~al.}(2020)Joshi, Mittal, and Pandey]{MR3914412}
Joshi, S.; Mittal, E.; Pandey, R.M.
\newblock On {E}uler type integrals involving extended {M}ittag-{L}effler functions.
\newblock {\em Bol. Soc. Parana. Mat. (3)} {\bf 2020}, {\em 38},~125--134.
	
\bibitem[Mainardi(2014)]{MAINARDI2}
Mainardi, F.
\newblock On some properties of the {M}ittag-{L}effler function
{$E_\alpha(-t^\alpha)$}, completely monotone for {$t>0$} with {$0<\alpha<1$}.
\newblock {\em Discrete Contin. Dyn. Syst. Ser. B} {\bf 2014}, {\em 19},~2267--2278.
\newblock
doi:{\changeurlcolor{black}\href{https://doi.org/10.3934/dcdsb.2014.19.2267}{\detokenize{10.3934/dcdsb.2014.19.2267}}}.
	
\bibitem[Schneider(1996)]{mettaglemma2}
Schneider, W.R.
\newblock Completely monotone generalized {M}ittag-{L}effler functions.
\newblock {\em Exposition. Math.} {\bf 1996}, {\em 14},~3--16.
	
\bibitem[Podlubny(1999)]{der1}
Podlubny, I.
\newblock {\em Fractional differential equations}; Academic Press, Inc., San
Diego, CA,  1999.
	
\bibitem[Sakamoto and Yamamoto(2011)]{solmittag}
Sakamoto, K.; Yamamoto, M.
\newblock Initial value/boundary value problems for fractional diffusion-wave
equations and applications to some inverse problems.
\newblock {\em J. Math. Anal. Appl.} {\bf 2011}, {\em 382},~426--447.
\newblock
doi:{\changeurlcolor{black}\href{https://doi.org/10.1016/j.jmaa.2011.04.058}{\detokenize{10.1016/j.jmaa.2011.04.058}}}.
	
\bibitem[Courant and Hilbert(1953)]{Aelliptic}
Courant, R.; Hilbert, D.
\newblock {\em Methods of mathematical physics. {V}ol. {I}}; Interscience
Publishers, Inc., New York, N.Y.,  1953; pp. xv+561.
	
\bibitem[Kato(1966)]{decomp}
Kato, T.
\newblock {\em Perturbation theory for linear operators}; Springer-Verlag New
York, Inc., New York,  1966.
	
\bibitem[Lazopoulos and Lazopoulos(2019)]{ref01}
Lazopoulos, K.A.; Lazopoulos, A.K.
\newblock On the Mathematical Formulation of Fractional Derivatives.
\newblock {\em Progr. Fract. Differ. Appl.} {\bf 2019}, {\em 5},~261--267.
\newblock
doi:{\changeurlcolor{black}\href{https://doi.org/10.18576/pfda/050402}{\detokenize{10.18576/pfda/050402}}}.
	
\end{thebibliography}
\end{document}